\documentclass[11pt]{article}
\usepackage[utf8]{inputenc}
\usepackage{authblk}

\usepackage{lmodern}
\usepackage[font=footnotesize]{caption}
\usepackage[font=footnotesize]{subcaption}
\usepackage[sf,bf,medium]{titlesec}

\usepackage[top=1in,bottom=1.2in,left=1.2in,right=1.2in]{geometry}

\usepackage{graphicx}
\usepackage{pgfplots}
\pgfplotsset{compat=1.14}

\usepackage{ulem}
\usepackage{mathtools}
\usepackage{booktabs}
\usepackage{amssymb}
\usepackage{amsmath}
\usepackage{color}
\usepackage{hyperref}

\usepackage[capitalize]{cleveref}
\crefname{equation}{Eq.}{Eqs.}
\Crefname{equation}{Equation}{Equations}
\crefrangelabelformat{equation}{(#3#1#4--#5#2#6)}
\crefmultiformat{equation}{Eqs. (#2#1#3}{, #2#1#3)}{#2#1#3}{#2#1#3}
\Crefmultiformat{equation}{Equations (#2#1#3}{, #2#1#3)}{#2#1#3}{#2#1#3}
\newcommand\pr[1]{\cref{#1}}

\newcommand{\Bext}{\mathbf{B}_{\mathrm{ext}}}
\newcommand{\Bplasma}{\mathbf{B}_{\mathrm{plasma}}}

\newcommand{\vct}{\mathbf}

\usepackage[sort,compress]{cite}

\numberwithin{equation}{section}

\newcommand{\lp}{\left(}
\newcommand{\rp}{\right)}

\newcommand{\D}[1]{\ensuremath{\,\mathrm{d}#1}}
\renewcommand{\vec}[1]{{\ensuremath{\mathbf{#1}}}}

\newcommand{\Boundary}{{\ensuremath{\Gamma}}}

\newcommand{\PAngle}{{\ensuremath{u}}}
\newcommand{\TAngle}{{\ensuremath{v}}}

\newcommand{\Normal}{{\vec{n}}}

\newcommand{\Nu}{\ensuremath{N_\PAngle}}                         
\newcommand{\Nv}{\ensuremath{N_\TAngle}}                         

\newcommand{\nexp}{{\ensuremath{\text{\sc{\small e-}}}}}         

\newcommand{\gmrestol}{\ensuremath{\epsilon_{\textsc{gmres}}}}   
\newcommand{\gmresiter}{\ensuremath{N_{\textsc{gmres}}}}         

\title{\bfseries \sffamily Efficient high-order singular quadrature
 schemes in magnetic fusion}

\author{Dhairya Malhotra\footnote{Email: \texttt{malhotra@cims.nyu.edu}},
  Antoine J. Cerfon\footnote{Email: \texttt{cerfon@cims.nyu.edu}},
  Michael O'Neil\footnote{Email: \texttt{oneil@cims.nyu.edu}},
  and Evan Toler\footnote{Email: \texttt{eht247@cims.nyu.edu}}
}
\date{\today}
\affil{Courant Institute, New York University\\
251 Mercer St, New York, NY 10012, USA}

\begin{document}

\maketitle
\begin{abstract}
  Several problems in magnetically confined fusion, such as the
  computation of exterior vacuum fields or the decomposition of the
  total magnetic field into separate contributions from the plasma and
  the external sources, are best formulated in terms of integral
  equation expressions. Based on Biot-Savart-like formulae, these
  integrals contain singular integrands. The regularization method
  commonly used to address the computation of various singular surface
  integrals along general toroidal surfaces is low-order accurate, and
  therefore requires a dense computational mesh in order to obtain
  sufficient accuracy. In this work, we present a fast, high-order
  quadrature scheme for the efficient computation of these integrals.
  Several numerical examples are provided demonstrating the
  computational efficiency and the high-order accurate convergence. A
  corresponding code for use in the community has been publicly
  released\footnote{\url{https://github.com/dmalhotra/BIEST}}.
\end{abstract}

\newpage

\section{Introduction}

Many magnetostatic calculations arising in magnetically confined
fusion applications can be efficiently expressed in terms of integral
equations~\cite{Freidberg1975,Lackner1976,Hirshman1986,Merkel1986,Merkel1987,Albanese1988,Chance1997,Atanasiu1999,Pustovitov2008,Lee2015,Marx2017}. A
standard example is the solution of the Neumann boundary value problem
corresponding to normal data on the magnetic field for Laplace's
equation in the interior or exterior of a toroidal
region~\cite{Chance1997,Freidberg1975,Freidberg1976,Gruber1981,Gruber19812,Merkel1982,Merkel1986,Pustovitov2008,Lazanja2011,Marx2017}. In
this case, the resulting integral equation for the vacuum field is
obtained by a direct application of Green's theorem. Another common
application of integral equations in this regime is
the virtual casing principle~\cite{stratton1941electromagnetic,Morozov1966,Shafranov1972,pustovitov2001magnetic,Hanson2015} used to
compute a
decomposition of the magnetic field into one field due to the plasma and
another field due to external sources
(e.g. coils)~\cite{Pustovitov2008,Lazanja2011,Drevlak2018}.

A well-known difficulty with integral formulations is that they often
involve integrals with singular kernels. Computing such integrals can
be challenging, both in terms of accuracy and in terms of
computational cost. The two challenges are intertwined, since
achieving high accuracy can require the use of computationally costly
quadrature schemes as well as fine meshes to properly resolve the
singularity. Several techniques have been proposed which address these
challenges for toroidally axisymmetric
boundaries~\cite{Chance2007,helsing2015explicit,oneil2018taylor}. However,
these advanced techniques do not apply to non-axisymmetric boundaries,
such as stellarator flux surfaces. In such situations, the standard
approach (at least in magnetic fusion codes)
is to remove the singular contribution in the integrand via
clever addition and subtraction of a known function with (to leading
order) the
appropriate singular behavior  which can be
integrated analytically. The regularized non-singular component is
then computed via the trapezoidal rule along the doubly periodic
non-axisymmetric
surface~\cite{Merkel1986,Lazanja2011,Drevlak2018}. This scheme
achieves the goal of removing numerical issues due to the singular
behavior of the kernel, but it is low-order accurate and thus generally obtains low precision and converges slowly. The reason for this is the following:
on a periodic interval, the trapezoidal rule
converges with order~$k+2$
for integrands which have $k$ continuous derivatives and
converges exponentially fast for integrands which are
analytic~\cite{Trefethen2014}.
The function resulting from the singularity
subtraction, which is integrated with the trapezoidal rule,
is bounded but not differentiable as it has a jump discontinuity.
This leads to slow convergence
of the trapezoidal rule. As a result of this slow convergence,
obtaining reasonable accuracy required for design and optimization
purposes requires a large number of quadrature nodes, in fact
larger than what is usually chosen for tractable
computations. This places a strong limitation on several solvers commonly
used in the stellarator community and which rely on this technique,
such as NESTOR~\cite{Merkel1986}, which is used in combination with
fixed-boundary stellarator equilibrium codes for the computation of
free-boundary equilibria~\cite{Hirshman1986,Strumberger1997}, as well as
both NESCOIL~\cite{Merkel1987} and REGCOIL~\cite{Landreman2017}.

While there exist several alternative quadrature rules based on
adaptive meshing and integration~\cite{Bremer_2012, Bremer_2013} and
analytic expansions/regularization~\cite{Kl_ckner_2013, Rachh_2017,
  Wala_2019, Rahimian_2017}, the geometries encountered in magnetic
fusion are most often doubly-periodic (e.g. flux surface in
stellarators) and therefore we wish to take advantage of the natural
global parameterization of these domains.  In this article, we present
a versatile, robust, and high-order accurate numerical method for
computing integrals with singular Green's functions (i.e. layer potentials)
along
general doubly-periodic toroidal surfaces.  Our method is also based
on singularity subtraction, but in contrast with the method described
above, we use a combination of a partition of unity and singularity
cancellation (via change of variables). The partition of unity
guarantees the smoothness of
the resulting non-singular integral which is computed using the
trapezoidal rule;  high-order accurate results are easily obtained.
Then, unlike the standard approach, the
remaining non-smooth integral cannot be evaluated analytically, but by
using a local polar coordinate system centered on the singularity it
can be computed to very high-order accuracy.  Our scheme converges
much faster than the standard singularity subtraction scheme of
Merkel~\cite{Merkel1987,Drevlak2018}, as we demonstrate numerically in
this article.

The structure of this article is as follows. In
Section~\ref{sec:two_problems} we review two applications of
magnetostatics commonly encountered in magnetic confinement fusion,
which involve the computation of integrals with singular kernels: the
separation of the magnetic field into contributions from the plasma
current and from external sources with the virtual casing principle,
and the computation of vacuum fields with Green's identity. These
applications allow us to illustrate our numerical scheme and test it,
but we insist that our algorithm is not limited to these situations.
In Section~\ref{sec:quadrature}, we discuss the low-order singularity
subtraction quadrature schemes currently used by the magnetic fusion
community, and present our high-order singularity subtraction
quadrature scheme. Then, in Section~\ref{sec:numerical_results} we
demonstrate the accuracy and speed of our high-order singularity
subtraction scheme via various numerical examples. Finally, we
summarize our results and point toward future areas of research in
Section~\ref{sec:summary}.

\section{Integral equations for magnetostatics in stellarators}
\label{sec:two_problems}

The high-order quadrature scheme we present in this article is
versatile and can be applied in a wide range of applications. However,
for the clarity of the presentation we focus here on the two most
common situations in stellarators in which our scheme may be used: (1)
the separation of the magnetic field into contributions from the
plasma current and from external sources (as needed, for example, in
coil design problems~\cite{Merkel1987,Landreman2017}), which we describe in Section \ref{subsec:normal-fields}, and (2) with
knowledge of the geometry of the plasma boundary, the
computation of vacuum magnetic fields~\cite{Hirshman1986,Merkel1986}, which we describe in Section \ref{subsec:vacuum-fields}.

\subsection{Calculating normal magnetic fields along a flux surface}
\label{subsec:normal-fields}
Let us denote by~$\Omega$ a toroidal region with smooth
boundary~$\Gamma$.  After solving fixed-boundary magnetohydrodynamic
(MHD) equilibrium problems, one has obtained the \textit{total}
magnetic field~$\mathbf{B}$ inside~$\Omega$.  Along~$\Gamma$, the
outermost flux surface, the equilibrium field is such
that~\mbox{$\mathbf{B} \cdot \mathbf{n} = 0$}, where $\mathbf{n}$ is
the outward unit normal vector along the boundary.  In general, the
total magnetic field~$\mathbf{B}$ is the superposition of two pieces:
one component is~$\Bplasma$, the magnetic field due to the plasma
current, and the other component is~$\Bext$, the magnetic field due to
the external coils.  In the plasma~$\mathbf{B} = \Bplasma + \Bext$,
and along the boundary~$\Gamma$ we have
that~$\Bext \cdot \mathbf{n} = - \Bplasma \cdot \mathbf{n}$ due to the
vanishing normal component of~$\mathbf{B}$.  In several applications,
such as coil design and the computation of the magnetic field in the
vacuum region outside the plasma, one is interested in the normal
component of~$\Bext$ along the plasma boundary.  To this end, using
the the Biot-Savart law we have that
\begin{equation}
  \mathbf{B}_{\mathrm{plasma}}(\mathbf{r})
  =\frac{\mu_{0}}{4\pi}\int_{\Omega}\frac{\mathbf{J}(\mathbf{r}')
    \times(\mathbf{r}-\mathbf{r}')}{|\mathbf{r}-\mathbf{r}'|^3} \,
  dv(\mathbf{r}'),
\label{eq:BS_volume}
\end{equation}
where~$\mathbf{J}$ is the current density in the plasma.
Eq.~(\ref{eq:BS_volume}) holds for any point $\mathbf{r}$ in $\Omega$
or on~$\Gamma$.  In particular, the above expression can be used to
compute~$\Bext \cdot\mathbf{n} = -\Bplasma \cdot\mathbf{n}$ on the
plasma boundary~$\Gamma$. However, there are two main computational
roadblocks when evaluating the above formula. First, because it is a
convolution, its numerical evaluation is a global calculation and
therefore typically slow. A much more
efficient strategy is to rely on the virtual casing
principle~\cite{Shafranov1972,Hanson2015} to rewrite this volume
integral as a surface integral. Using a Green's-like identity for the
magnetic field, the normal component of the
magnetic field due to the plasma can be expressed as:
\begin{equation}
  \mathbf{B}_{\mathrm{plasma}}(\mathbf{r})\cdot\mathbf{n}(\mathbf{r})
  =
  \frac{1}{4\pi} \, \mathbf{n}(\mathbf{r}) \cdot
  \int_{\Gamma}\frac{(\mathbf{n}(\mathbf{r}')\times\mathbf{B}(\mathbf{r}'))\times(\mathbf{r}-\mathbf{r}')}{|\mathbf{r}-\mathbf{r}'|^3}
  \, da(\mathbf{r}') ,
\label{eq:BS_surface}
\end{equation}
where~$\mathbf{r}$ is an observation point on~$\Gamma$.
The second computational roadblock, which applies to both the volume formulation in Eq. (\ref{eq:BS_volume}) and the surface formulation given by Eq.~(\ref{eq:BS_surface}), is that the integrand is singular when the
point~$\mathbf{r}'$ is the same as the observation point
$\mathbf{r}$. This difficulty is precisely the one we address in this article.

Some numerical codes choose to evaluate the quantity in
Eq.(\ref{eq:BS_surface}) in two steps~\cite{Drevlak2018}. First,
the vector potential~$\mathbf{A}_{\mathrm{plasma}}$ is computed by
evaluating:
\begin{equation}
  \mathbf{A}_{\mathrm{plasma}}(\mathbf{r})=\frac{1}{4\pi}\int_{\Gamma}
  \frac{\mathbf{n}(\mathbf{r}')\times\mathbf{B}(\mathbf{r}')}
  {|\mathbf{r}-\mathbf{r}'|} \, da(\mathbf{r}').
\label{eq:vecpot_surface}
\end{equation}
Then, subsequent Fourier differentiation is used to compute the
tangential derivatives of the tangential components
of~$\mathbf{A}_{\mathrm{plasma}}$, from
which~$\mathbf{B}_{\mathrm{plasma}}\cdot\mathbf{n}$ can be
calculated. Similarly,
Eq.~(\ref{eq:vecpot_surface}) also involves a singular
integrand when the evaluation point~$\mathbf{r}'$ is the same as the
observation point~$\mathbf{r}$. The kernel appearing in
Eq.~(\ref{eq:vecpot_surface}) is the kernel appearing in Laplace
single-layer potentials, while the kernel appearing in
Eq.~(\ref{eq:BS_surface}) is closely related to the kernel for
double-layer potentials~\cite{guentherlee1995}.  The
method we propose in this article applies to computing the integrals
in
both Eq.~(\ref{eq:BS_surface}) and Eq.~(\ref{eq:vecpot_surface}).

\subsection{Computing vacuum magnetic fields}
\label{subsec:vacuum-fields}

Consider a magnetic field $\mathbf{B}_{0}$ satisfying
\begin{equation}
  \nabla\times\mathbf{B}_{0}=\mathbf{J},
  \qquad\nabla\cdot\mathbf{B}_{0}=0,
  \label{eq:B0Vac}
\end{equation}
in $\mathbb R^3 = \Omega\cup\Omega^{\mathrm{c}}$, where as
before~$\Omega$ is the plasma domain,~$\Omega^{\mathrm{c}}$ is the
exterior of~$\Omega$, and~$\mathbf{J}$ is the total
current density (both in the plasma and in~$\Omega^{\mathrm{c}}$,
i.e. in coils exterior to the plasma). Along a general  boundary~$\Gamma=\partial\Omega$,
the normal component of~$\mathbf{B}_0$ may be non-zero. Additionally,
however, we can assume there exists a vacuum field~$\mathbf{B}_{\mathrm{vac}}$
in~$\Omega^c$ with zero poloidal circulation, such that the total field 
$\mathbf{B}_{\mathrm{tot}} = \mathbf{B}_{0}+\mathbf{B}_{\mathrm{vac}}$ 
satisfies the following equations \cite{Hirshman1986,Merkel1986}
\begin{equation}
  \begin{aligned}
    \nabla\times\mathbf{B}_{\mathrm{tot}}&=\mathbf{J},
    &\qquad &\text{in }\Omega^{c}, \\
    \nabla\cdot\mathbf{B}_{\mathrm{tot}} &=0, & &\text{in }\Omega^{c},
    \\
    \mathbf{n}\cdot\mathbf{B}_{\mathrm{tot}} &= 0, &
    &\text{on }\Gamma=\partial\Omega.
  \end{aligned}
\end{equation}
Physically, the field $\mathbf{B}_{\mathrm{vac}}$ can be viewed as being due to surface currents flowing on $\Gamma$, which are precisely such that the normal component of $\mathbf{B}_{\mathrm{tot}}$ vanishes on $\Gamma$. These surface currents are not current sources included in $\mathbf{J}$ in (\ref{eq:B0Vac}), and are such that they do not contribute to a net toroidal current. This is what we mean by ``zero poloidal circulation''. Because $\mathbf{B}_{\mathrm{vac}}$ has zero poloidal circulation, it can be uniquely written
as~$\mathbf{B}_{\mathrm{vac}}=\nabla\Phi$, where $\Phi$ is a
single-valued potential which satisfies the following exterior Laplace equation with Neumann boundary conditions~\cite{Hirshman1986,Merkel1986}:
\begin{equation}
  \begin{aligned}
    \Delta\Phi &=0, &\qquad &\text{in }\Omega^{c}, \\
    \mathbf{n}\cdot\nabla\Phi &= -\mathbf{n}\cdot\mathbf{B}_{0}, &
    &\text{on }\Gamma, \\
    \Phi &\to 0, & &\text{as } r \to \infty.
  \end{aligned}
  \label{eq:neumannproblem}
\end{equation}

As discussed in the introduction, this boundary value problem can be
rewritten in integral form via a direct application of Green's
identity.  The potential~$\Phi$ can be shown to satisfy the following
equation along~$\Gamma$:
\begin{equation}
  \Phi(\mathbf{r}) -
  \frac{1}{2\pi}\int_{\Gamma}\mathbf{n}(\mathbf{r}') \cdot
  \nabla_{\mathbf{r}'}
  \left(\frac{1}{|\mathbf{r}-\mathbf{r}'|}\right) \,
  \Phi(\mathbf{r}') \, da(\mathbf{r}') =
  -\frac{1}{2\pi}\int_{\Gamma}
  \frac{\mathbf{n}(\mathbf{r}') \cdot
    \nabla\Phi(\mathbf{r}')}{|\mathbf{r}-\mathbf{r}'|} \, da(\mathbf{r}').
\label{eq:greensbasic}
\end{equation}
With the Neumann boundary condition given
in~\eqref{eq:neumannproblem}, the right-hand side
of~Eq.(\ref{eq:greensbasic}) can be evaluated directly so
that~(\ref{eq:greensbasic}) is an integral equation for~$\Phi$
on~$\Gamma$. Furthermore, once $\Phi$ has been computed on~$\Gamma$,
it can be evaluated anywhere in $\Omega^{c}\setminus\Gamma$ using
virtually the same Green's identity~\cite{Merkel1986}, except
evaluated off-surface:
\begin{equation}
  \Phi(\mathbf{r})=\frac{1}{4\pi}\int_{\Gamma}
  \mathbf{n}(\mathbf{r}')\cdot\nabla_{\mathbf{r}'}
  \left(\frac{1}{|\mathbf{r}-\mathbf{r}'|}\right)\,
  \Phi(\mathbf{r}') \, da(\mathbf{r}')-\frac{1}{4\pi}\int_{\Gamma}
  \frac{1}{|\mathbf{r}-\mathbf{r}'|}\mathbf{n}(\mathbf{r}')\cdot
  \nabla\Phi(\mathbf{r}') \, da(\mathbf{r}').
\label{eq:greensoutside}    
\end{equation}
Free boundary MHD equilibrium computations often rely on the
formulation we just presented~\cite{Hirshman1986}. A nearly identical
integral equation approach can also be used for the computation of
interior vacuum fields~\cite{Merkel1986}. We observe that as in the
previous subsection, numerical schemes based on such integral equation
formulations need to accurately and efficiently treat the singular
integrands appearing in the integrals. Here, two singular kernels must
be addressed: the singularity of the single-layer potential type, with
$1/|\mathbf{r}-\mathbf{r}'|$ as the singular kernel, and a
singularity of the double-layer type, i.e. the normal derivative
 of~$1/|\mathbf{r}-\mathbf{r}'|$ along the surface.

 At this point, we emphasize that the singularity subtraction scheme
 currently used in the magnetic fusion community and the new
 high-order singularity subtraction scheme we present in this article
 only apply to \textit{on-surface evaluations}, as needed in
 Eq.(\ref{eq:greensbasic}). To the best of our knowledge, methods for
 \textit{off-surface evaluations}, as required in
 Eq.(\ref{eq:greensoutside}), have not been proposed in the magnetic
 fusion community. In fact, designing efficient and accurate numerical
 methods for \textit{off-surface evaluations} near arbitrary surfaces
 in three dimensions remains an open problem in applied mathematics,
 although the pioneering work of~\cite{Zhao_2010, Corona_2017,
   Tlupova_2013, Kl_ckner_2013, Wala_2019, Af_Klinteberg_2014,
   Siegel_2018, Ying_2006, Rahimian_2017, Carvalho_2018a,
   Carvalho_2018b} must be highlighted. Fortunately, for many
 applications, including free-boundary MHD equilibrium calculations,
 the vacuum magnetic field is only needed \textit{on the magnetic
   surface}, so that only Eq.(\ref{eq:greensbasic}) needs to be
 solved. We will focus on \textit{on-surface evaluations} in this
 article.

\section{Quadrature for surface integrals}
\label{sec:quadrature}

\subsection{Singularity subtraction}
\label{sec:loworder}

In this section, we review the method originally proposed by Merkel to
address the challenges associated with the singularity appearing in
the integrands of the integral formulations presented above. This
scheme is a workhorse of several currently widely used MHD
codes~\cite{Hirshman1986,Merkel1986,Drevlak2018}. The following
presentation is brief since it has already been discussed in detail by
Merkel~\cite{Merkel1986}, by Lazanja~\cite{Lazanja2011}, and by
Drevlak \textit{et al.}~\cite{Drevlak2018}.

In the simplest case of the layer potentials of the previous section,
one must numerically compute an integral of the form:
\begin{equation}
    \Phi(\vct{r}) = \frac{1}{4\pi} 
    \int_\Gamma \frac{\sigma(\vct{r}')}{|\vct{r} - \vct{r}'|} \, 
    da(\vct{r}'),
\end{equation}
where~$\sigma$ is some smooth function defined over the doubly-periodic
surface~$\Gamma$. If the surface~$\Gamma$ is parameterized by the
variables~$u,v$, i.e. if~$\Gamma: [0,2\pi) \times [0,2\pi) \to \mathbb R^3$,
then for targets~\mbox{$\vct{r} = \vct{r}(u,v) \in \Gamma$}, 
the above integral can be written:
\begin{equation}
    \Phi(u,v) = \frac{1}{4\pi} \int_0^{2\pi} \int_0^{2\pi}
    \frac{\sigma(u',v')}{|\vct{r}(u,v) - \vct{r}(u',v') |} \, 
    g(u',v') \, du' \, dv',
\end{equation}
where~$g$ is the metric tensor along~$\Gamma$ and in a small abuse of 
notation,~$\vct{r}$ denotes the parameterization of~$\Gamma$.
Take note that the above integral is~\textit{not} exactly a convolution
in the~$u,v$ variables, as it
is an integral along a curved surface.
Taylor expanding~$\vct{r}'$ about the point~$(u,v)$ we have that
\begin{equation}
    \vct{r}(u',v') \approx \vct{r} + \vct{r}_u \, (u'-u) + 
    \vct{r}_v \, (v'-v)
\end{equation}
where
\begin{equation}
    \vct{r}_u = \frac{\partial \vct{r}}{\partial u}(u,v), \qquad
    \vct{r}_v = \frac{\partial \vct{r}}{\partial v}(u,v).
\end{equation}
Therefore, to leading order
\begin{equation}
  \begin{aligned}
    | \vct{r}(u,v) - \vct{r}(u',v') | &\approx 
    \sqrt{\vct{r}_u \cdot \vct{r}_u \, (u'-u)^2 + 
      2 \vct{r}_u \cdot \vct{r}_v \, (u'-u)(v'-v) + 
      \vct{r}_v \cdot \vct{r}_v \, (v'-v)^2 } \\
    &= R(u,v,u',v').
  \end{aligned}
\end{equation}
Next, adding and subtracting we have that
\begin{multline}
  \Phi(u,v) = \frac{1}{4\pi} \int\int \lp 
  \frac{\sigma(u',v') \, g(u',v')}{|\vct{r}(u,v) - \vct{r}(u',v') |}
  - \frac{\sigma(u,v) \, g(u,v) }{R(u,v,u',v')} \rp 
  du' \, dv' + \\
  \frac{1}{4\pi} \int\int  
  \frac{\sigma(u,v) \, g(u,v) }{R(u,v,u',v')}
  du' \, dv'.
\end{multline}
In this form, the singularity has been subtracted from the integrand
in the first integral, which can therefore be computed using the
trapezoidal rule. The expected order of
convergence of the trapezoidal rule for this integral is only
2nd-order, however, despite the doubly periodic nature of the
integrand. Indeed, while the leading order singularity has been
subtracted,
the integrand has a jump discontinuity at a point
and therefore exponential convergence of the trapezoidal rule is not
expected, and not observed in practice.
With trapezoidal rule, exponential converge is observed only for
periodic integrands which are also analytic \cite{Trefethen2014}.
As shown in
\cite{Drevlak2018}, the second integral above can be analytically
integrated in one of the coordinates. The resulting one-dimensional
integral still has a singular integrand, which can be treated with the
same technique: addition and subtraction of a term with the
appropriate singular behavior, and analytic integration of the
additional integral introduced. We omit the actual formulas, as they
are rather complicated, but refer the reader to~\cite{Drevlak2018} for
the specific case of Biot-Savart-like kernels.

The brief description we just gave follows the approach of Drevlak
\textit{et al.} quite closely. We chose it as an illustration of the
standard singularity subtraction scheme because it captures the
central idea of the scheme with the minimal amount of notation and
relatively few cumbersome expressions. We note that the original
singularity subtraction scheme proposed in \cite{Merkel1986}
is very similar in spirit, and has the same low-order convergence
properties, but differs in two ways. First, in \cite{Merkel1986}
the goal is to compute the vacuum magnetic field using a spectral
scheme based on Fourier series to solve Eq.(\ref{eq:greensbasic}). One
therefore needs to compute the Fourier transforms of the integrals in
Eq.(\ref{eq:greensbasic}), instead of the integrals
themselves. Second, to compute these Fourier transforms, the author does
not add and subtract the term $\sigma(u,v) \, g(u,v) /R(u,v,u',v')$
for the single-layer potential integral but instead a well-chosen
expression which has the same singular behavior as that term, and for
which he can remarkably compute the Fourier transform analytically. In
the convergence tests in Section \ref{ss:result-dbl-layer}, in which
we compare our high-order scheme with the low-order singularity
subtraction technique, our implementation of the analytic singularity
subtraction scheme is based on the approach presented in~\cite{Merkel1986}. We have not
implemented the approach by Drevlak \textit{et al.}~\cite{Drevlak2018} but we expect
similar performance of that scheme, based on the mathematical reasons
we gave previously.

Finally, we observe that we focused on kernels of the form appearing
in Laplace single-layer potentials in this subsection. We did so
merely for the simplicity and conciseness of the
presentation. Analytic singularity subtraction schemes also exist for kernels of the form of double-layer potentials, as described
explicitly and used in \cite{Merkel1986}. Our convergence tests in
Section~\ref{ss:result-dbl-layer} will precisely focus on  kernels
similar to those appearing in the double-layer potential.

\subsection{A high-order singularity subtraction scheme}
\label{ss:high-order-quad}

The previous singularity subtraction scheme is relatively simple to
implement due to the existence of analytic formulae for the integral
of the subtracted part. However, since only the leading order term in
the singularity is canceled via addition/subtraction,
only low-order convergence is observed,
and the actual accuracy obtained is modest at best. In this section,
we present a high-order singularity subtraction scheme which addresses
these limitations, and can be used to increase the speed and accuracy
of the fusion codes which rely on the numerical evaluation of
integrals of the forms discussed above; variants of the scheme were
originally used in \cite{Bruno_2001a,Bruno_2001b,Ying_2006} and an
improved version was used in~\cite{Malhotra2019} as well as for the
numerical results in the present article.

\newcommand{\KernelFn}{K}
\newcommand{\MetricTensor}{\vec{G}}
\newcommand{\pou}{\ensuremath{\chi}}
\newcommand{\SurfCoord}{\vec{r}}
\newcommand{\POUSupport}{\eta}
\newcommand{\DiffAreaElem}[2]{g(#1,#2)}

\begin{figure}
  \resizebox{\textwidth}{!}{%
  \begin{tikzpicture}[scale=1]
    \node[inner sep=0pt] at (0.0, 0.0) {\includegraphics[width=0.27\textwidth]{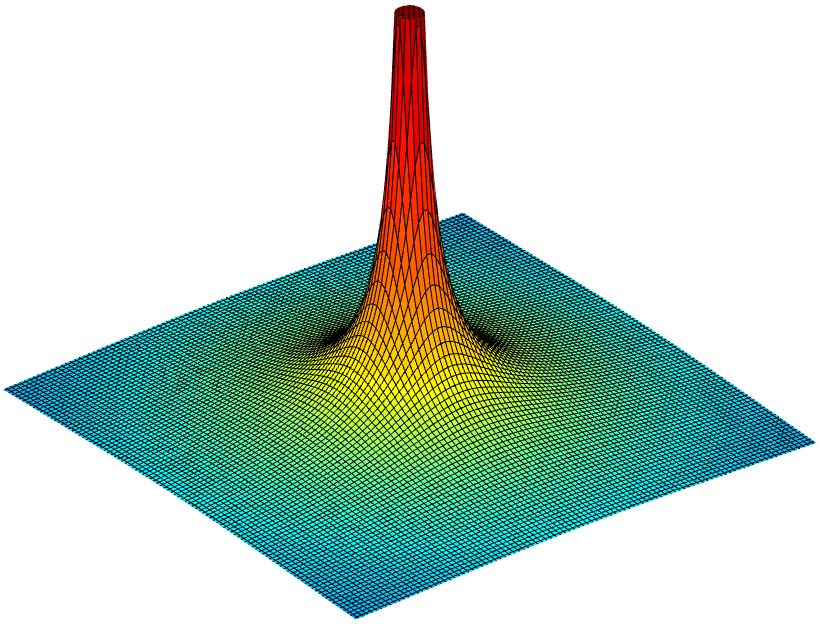}};
    \node[inner sep=0pt] at (3.8,-1.5) {\includegraphics[width=0.27\textwidth]{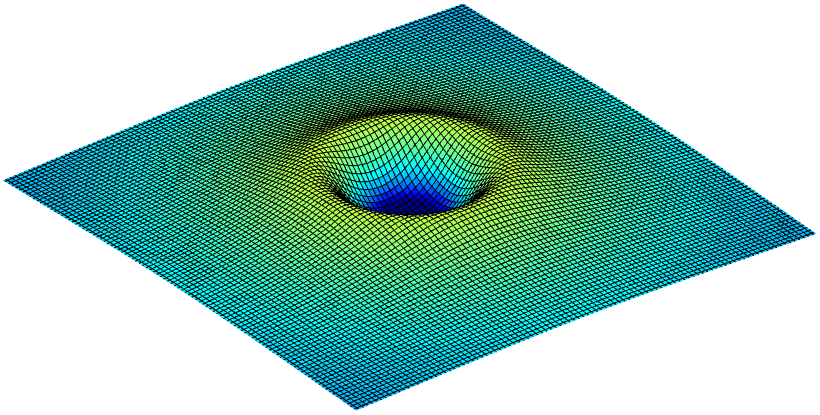}};
    \node[inner sep=0pt] at (3.1, 0.7) {\includegraphics[width=0.08775\textwidth]{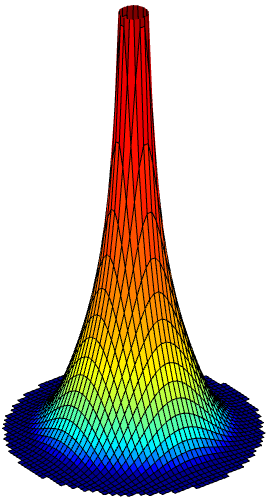}};
    \node[inner sep=0pt] at (7.4,-0.3) {\includegraphics[width=0.23\textwidth]{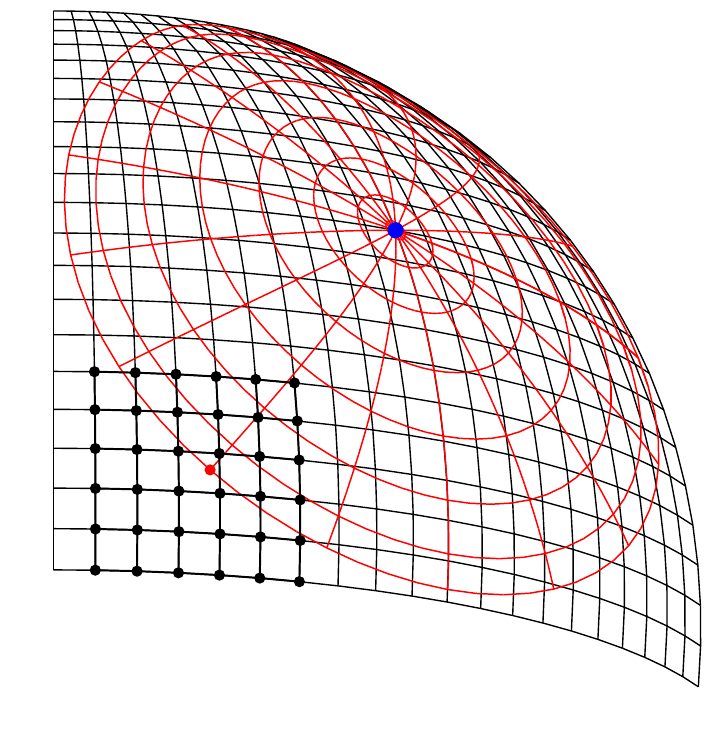}};
  \end{tikzpicture}}
  \caption{\label{f:high-order-quad}
  The Laplace single-layer kernel function is shown on the left.
  We use a partition of unity to split it into a $C^\infty$ smooth part (bottom) and a local singular part (top).
  The convolution with the smooth kernel function is computed using trapezoidal quadrature, leading to exponential convergence.
  The integral containing the singular part of the kernel is computed in polar coordinates after interpolating from the regular grid to a polar grid using local polynomial interpolation (right). The coordinate transformation cancels the kernel singularity, and trapezoidal quadrature for the angular direction and Gauss-Legendre quadrature for the radial direction then also leads to spectral convergence.
  }
\end{figure}

To this end we consider the problem of computing singular integrals
over a double-periodic surface $\Boundary$,
\begin{align}
  \Phi(\SurfCoord(\PAngle,\TAngle)) &= \int_\Boundary \KernelFn(\SurfCoord(\PAngle,\TAngle) - \SurfCoord') ~ \sigma(\SurfCoord') ~ da(\SurfCoord') \nonumber \\
                                    &= \int_{0}^{2\pi} \int_{0}^{2\pi} \KernelFn(\SurfCoord(\PAngle,\TAngle) - \SurfCoord(\PAngle',\TAngle')) ~ \sigma(\PAngle',\TAngle') ~ \DiffAreaElem{\PAngle'}{\TAngle'} ~ d\PAngle' ~ d\TAngle'
  \label{e:generic-boundary-integ}
\end{align}
where $g$, $\sigma$, and $\SurfCoord$ are as defined previously in
\pr{sec:loworder} and $\KernelFn$ is a singular kernel function such
as the single-layer Laplace kernel
function~($\KernelFn(\vec{r}) = 1 / 4\pi |\vec{r}|$) or
its derivatives.  As previously mentioned, such integrals are
challenging to evaluate due to the singularity in the kernel function.
The high-order singularity cancellation scheme we propose uses a
partition of unity to split the singular integral in
\pr{e:generic-boundary-integ} into an integral with $C^\infty$
integrands (as opposed to only bounded integrands in the case of the
low-order scheme discussed \pr{sec:loworder}) and a local singular
integral,
\begin{multline}
  \Phi(\SurfCoord(\PAngle,\TAngle)) = \int_{0}^{2\pi} \int_{0}^{2\pi} \KernelFn_{\text{smooth}}(\PAngle',\TAngle') ~ \sigma(\PAngle',\TAngle') ~ \DiffAreaElem{\PAngle'}{\TAngle'} ~ d\PAngle' ~ d\TAngle' + \\
                                  \int_{-\POUSupport}^{\POUSupport} \int_{-\POUSupport}^{\POUSupport} \KernelFn_{\text{singular}}(\PAngle',\TAngle') ~ \sigma(\PAngle',\TAngle') ~ \DiffAreaElem{\PAngle'}{\TAngle'} ~ d\PAngle' ~ d\TAngle'
  \label{e:high-order-split}
\end{multline}
where
\begin{align}
  \KernelFn_{\text{singular}}(\PAngle',\TAngle') &= \KernelFn(\SurfCoord(\PAngle,\TAngle) - \SurfCoord(\PAngle',\TAngle')) ~ \pou\!\left(\frac{1}{\POUSupport}\sqrt{(\PAngle-\PAngle')^2+(\TAngle-\TAngle')^2}\right) \\
  \KernelFn_{\text{smooth}}(\PAngle',\TAngle')   &= \KernelFn(\SurfCoord(\PAngle,\TAngle) - \SurfCoord(\PAngle',\TAngle')) ~ \left( 1 - \pou\!\left(\frac{1}{\POUSupport}\sqrt{(\PAngle-\PAngle')^2+(\TAngle-\TAngle')^2}\right) \right).
\end{align}
Here, $\pou$ can be any smooth function such that $\pou(x) = 1$ in a neighborhood around zero with compact support on [-1,1].
This ensures that $\KernelFn_{\text{smooth}}$ is smooth and $\KernelFn_{\text{singular}}(\PAngle',\TAngle')$ has compact support on $(\PAngle',\TAngle') \in [-\POUSupport,\POUSupport]^2$.
The first integral in \pr{e:high-order-split} has smooth integrands
and can be evaluated using trapezoidal quadrature which gives
exponential convergence.
Unlike the low order scheme, the singular integral is not known
analytically; however, it can be computed numerically to high-order
accuracy by computing the integral in polar coordinates,
\begin{multline}
  \int_{-\POUSupport}^{\POUSupport} \int_{-\POUSupport}^{\POUSupport}
  \KernelFn_{\text{singular}}(\PAngle',\TAngle') ~
  \sigma(\PAngle',\TAngle') ~ \DiffAreaElem{\PAngle'}{\TAngle'} ~
  d\PAngle' ~ d\TAngle' \\
  =
  \int_{0}^{\POUSupport} \int_{0}^{2\pi} \rho
  \KernelFn_{\text{singular}}(\rho,\theta) ~ \sigma(\rho,\theta) ~
  \DiffAreaElem{\rho}{\theta} ~ d\theta ~ d\rho ,
\end{multline}
where we have applied the change of variables
$\PAngle = \rho \cos \theta$ and $\TAngle = \rho \sin \theta$.
The Jacobian of the coordinate transformation cancels kernel
singularities
of the form $1/(\mathbf{r}-\mathbf{r}')$
and the integral can then be computed using trapezoidal
quadrature in the angular direction $\theta \in [0,2\pi)$ and
Gauss-Legendre quadrature in the radial direction
$\rho \in [0,\POUSupport]$.
In \cite{Ying_2006} a proof was given to show that this scheme also works for the double-layer Stokes pressure kernel and in general for other integrals which must be understood in principle value sense, such as the potential due to the Biot-Savart kernel shown in \pr{eq:BS_surface}.
The two central steps of our scheme,
namely the split of the singular kernel into a smooth kernel and a
singular kernel with compact support, and the coordinate
transformation and interpolation to a polar grid, are shown in Figure
\ref{f:high-order-quad}.
The details of the complete algorithm and the
performance optimizations can be found in \cite{Malhotra2019}.
The optimal choice for size of the support $\POUSupport$ (and also the quadrature order for integrating in polar coordinates) can be guessed from the number of grid points required to resolve $\pou$ (shown in Fig.~4 in \cite{Malhotra2019}). The scheme is spectrally accurate in each parameter and the overall error is the maximum of the error from each component (i.e. the surface and density discretization which depends on the grid resolution; and the accuracy of the smooth and singular quadratures which depend on $\POUSupport$ and the order of the quadrature in polar coordinates).
The numerical results in Table~3 of \cite{Malhotra2019} show the accuracy of the scheme for different choices of these parameters. The cost of the scheme as a function of the different parameters is given in Table~2 in \cite{Malhotra2019}.

\section{Numerical results}
\label{sec:numerical_results}
In this section we present several numerical results demonstrating the
high-order accurate convergence of the quadrature scheme presented in
Section~\ref{ss:high-order-quad}.
In \pr{ss:result-dbl-layer}, we compare the high-order and the low-order quadrature schemes using a simple Green's identity test.
In \pr{ss:result-ext-field} we show numerical results for computing the normal component of the external field in a non-axisymmetric geometry.
Finally, in \pr{subsec:results-vacuum-field} we present results for computing vacuum magnetic fields.

\subsection{Double-layer potential on an arbitrary toroidal surface \label{ss:result-dbl-layer}}
As we have seen, from Green's identity, a harmonic function $\Phi(\vec{r})$ in $\Omega$ can be represented by the sum of a single-layer and double-layer potential,
\begin{equation}
  \Phi(\vec{r}) = -\frac{1}{2\pi} \int_{\Gamma} \Phi(\vec{r'}) \frac{\vec{n}(\vec{r'}) \cdot \left( \vec{r}-\vec{r'} \right)} {|\vec{r}-\vec{r'}|^3} da(\mathbf{r}')
                  + \frac{1}{2\pi} \int_{\Gamma} \frac{\partial \Phi}{\partial \vec{n}}(\vec{r'}) \frac{1} {|\vec{r}-\vec{r'}|} da(\mathbf{r}')
                  ,
\end{equation}
where $\vec{r}$ is a point on $\Boundary = \partial\Omega$.
For the case $\Phi = 1$ in $\Omega$, this results in the identity,
\begin{equation}
  \begin{aligned}
  e(\vec{r}) &=  1 + \frac{1}{2\pi} \int_{\Gamma}
  \frac{\vec{n}(\vec{r'}) \cdot \left( \vec{r}-\vec{r'} \right)}
  {|\vec{r}-\vec{r'}|^3} da(\mathbf{r}') \\
  &= 0.
  \end{aligned}
\end{equation}
We numerically evaluate the above layer potentials and plot the error
$\|e\|_\infty$ versus mesh refinement in \pr{f:dbl-conv}.  We show
results for two different surface geometries: first for an
axisymmetric torus with a circular cross section and second for a
non-axisymmetric surface which approximates the shape of the plasma
boundary in the W7-X stellarator.  In both cases we compare the
low-order singularity subtraction scheme of Merkel~\cite{Merkel1986}
and the high-order partition of unity based scheme. In the low-order
integration scheme with staggering, the evaluation points are at
integer points on the surface, and the quadrature nodes are at
half-integer points \cite{Drevlak2018}; in contrast, in the low-order
scheme without staggering, evaluation points and quadrature points are
the same. We observe that the low-order scheme converges with a
second-order rate of convergence while the high-order scheme has a
numerical convergence rate of greater than 10th order.
For exponential convergence, the slope of the $\log-\log$ plot
of the error becomes steeper with increasing grid resolution; therefore,
we expect the numerically observed order to continue increasing at higher
accuracies and the observed results are consistent with exponential convergence.
We also
note that the W7-X geometry has regions of high curvature and
therefore requires finer surface discretizations to achieve the same
quadrature accuracy as the axisymmetric geometry; however, the
convergence rates in both cases are as expected from the theory.

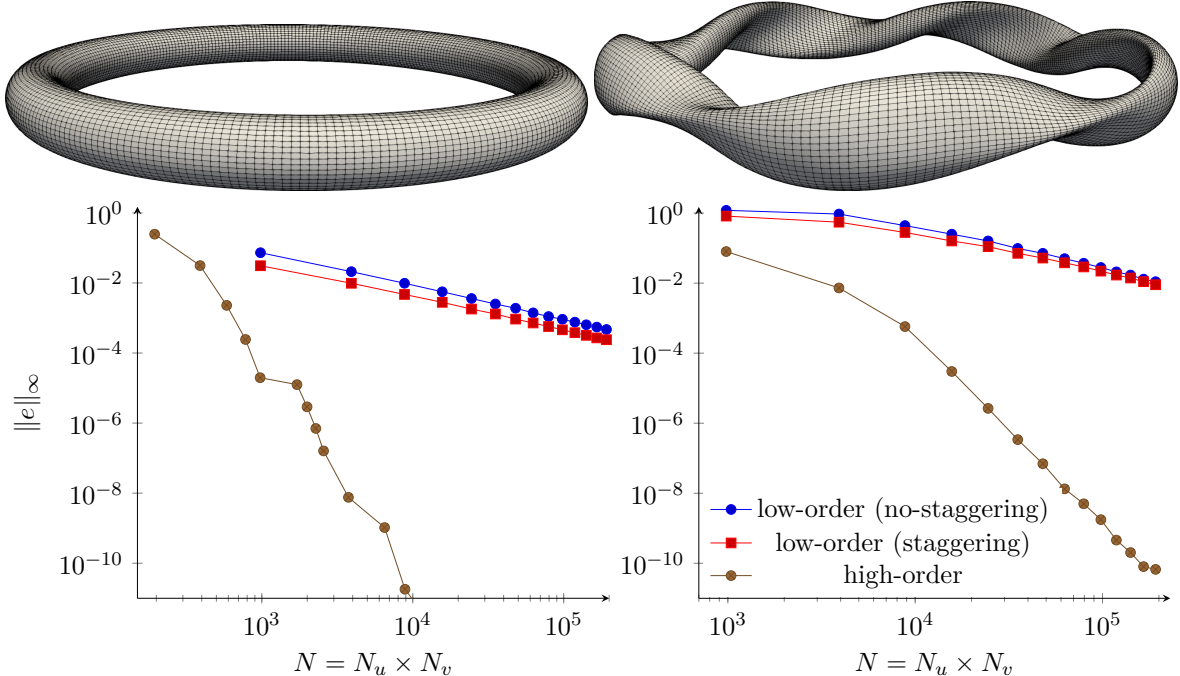
\begin{figure}[t]
  \includegraphics[width=0.5\textwidth]{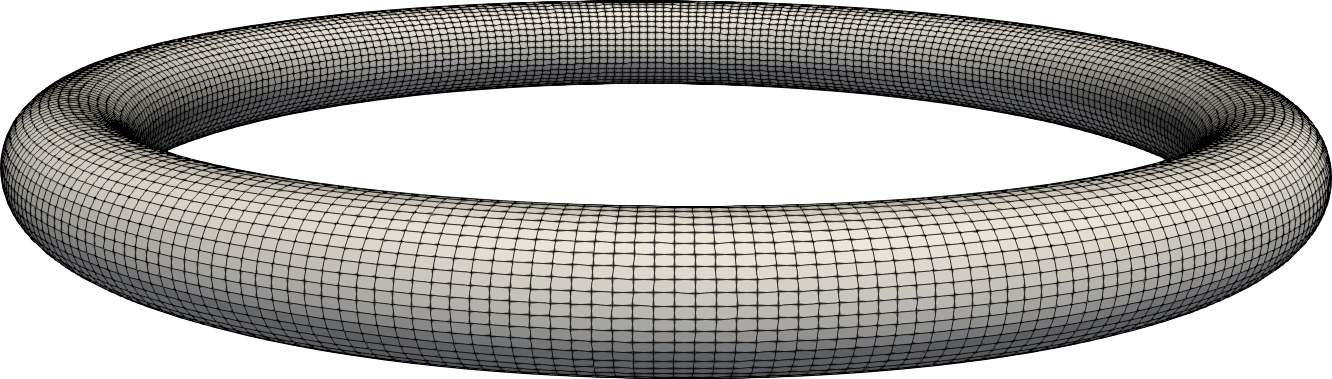}
  \includegraphics[width=0.5\textwidth]{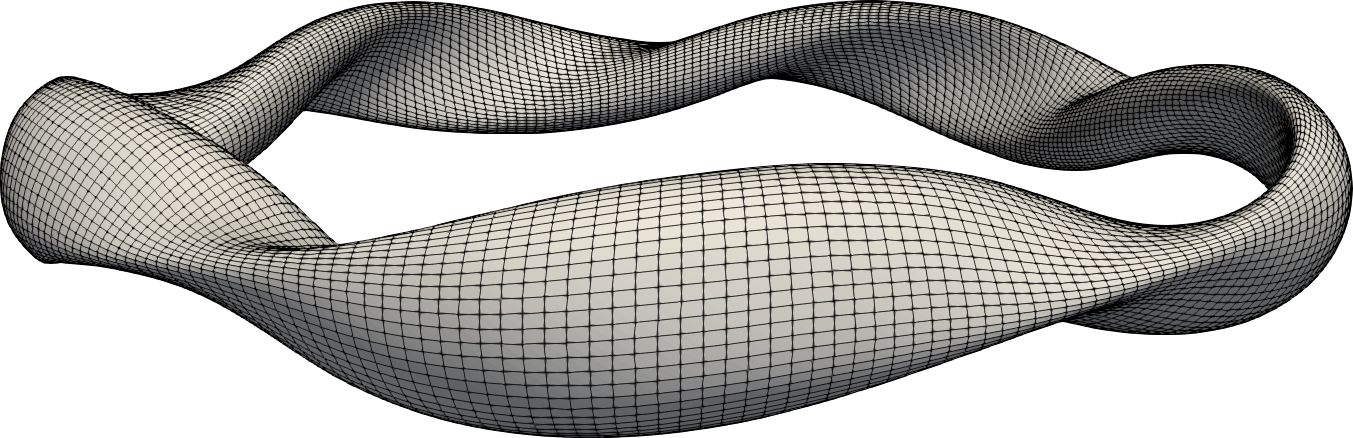}
  \resizebox{0.52\textwidth}{!}{\begin{tikzpicture}
    \begin{axis}[
        xmode=log,
        ymode=log,
        axis lines=left,
        xlabel={$N = N_u \times N_v$},
        ylabel={$\|e\|_\infty$},
        ymin=1e-11,
        ymax=1.5,
        xmin=1.5e2,
        xmax=2e5,
        legend style={draw=none},
      ]

      \addplot coordinates {
        (   980, 7.3e-02)
        (  3920, 2.1e-02)
        (  8820, 9.8e-03)
        ( 15680, 5.6e-03)
        ( 24500, 3.6e-03)
        ( 35280, 2.5e-03)
        ( 48020, 1.9e-03)
        ( 62720, 1.4e-03)
        ( 79380, 1.1e-03)
        ( 98000, 9.2e-04)
        (118580, 7.6e-04)
        (141120, 6.4e-04)
        (165620, 5.5e-04)
        (192080, 4.7e-04)
      };

      \addplot coordinates {
         (   980, 3.1e-02)
         (  3920, 9.8e-03)
         (  8820, 4.7e-03)
         ( 15680, 2.8e-03)
         ( 24500, 1.8e-03)
         ( 35280, 1.3e-03)
         ( 48020, 9.3e-04)
         ( 62720, 7.2e-04)
         ( 79380, 5.7e-04)
         ( 98000, 4.6e-04)
         (118580, 3.8e-04)
         (141120, 3.2e-04)
         (165620, 2.7e-04)
         (192080, 2.4e-04)
      };

      \addplot coordinates {
        (  195, 0.246547   )
        (  390, 0.0313301  )
        (  585, 0.0022925  )
        (  780, 0.00024359 )
        (  975, 1.97686e-05)
        ( 1710, 1.2541e-05 )
        ( 1995, 2.91152e-06)
        ( 2280, 7.03463e-07)
        ( 2565, 1.60672e-07)
        ( 3750, 7.59613e-09)
        ( 6510, 1.04357e-09)
        ( 8880, 1.80728e-11)
        (12900, 2.14129e-12)
      };


    \end{axis}
  \end{tikzpicture}}
  \resizebox{0.48\textwidth}{!}{\begin{tikzpicture}
    \begin{axis}[
        xmode=log,
        ymode=log,
        axis lines=left,
        xlabel={$N = N_u \times N_v$},
        ymin=1e-11,
        ymax=1.5,
        xmin=7e2,
        xmax=2.3e5,
        legend style={draw=none,at={(axis cs:730,1.1e-11)},anchor=south west},
      ]

      \addplot coordinates {
        (   980, 1.2e+00)
        (  3920, 9.4e-01)
        (  8820, 4.4e-01)
        ( 15680, 2.5e-01)
        ( 24500, 1.6e-01)
        ( 35280, 9.8e-02)
        ( 48020, 7.1e-02)
        ( 62720, 5.0e-02)
        ( 79380, 3.7e-02)
        ( 98000, 2.8e-02)
        (118580, 2.1e-02)
        (141120, 1.7e-02)
        (165620, 1.3e-02)
        (192080, 1.1e-02)
      };
      \addlegendentry{low-order (no-staggering)}

      \addplot coordinates {
        (   980, 8.2e-01)
        (  3920, 5.5e-01)
        (  8820, 2.8e-01)
        ( 15680, 1.6e-01)
        ( 24500, 1.1e-01)
        ( 35280, 7.1e-02)
        ( 48020, 5.2e-02)
        ( 62720, 3.8e-02)
        ( 79380, 2.9e-02)
        ( 98000, 2.2e-02)
        (118580, 1.7e-02)
        (141120, 1.4e-02)
        (165620, 1.1e-02)
        (192080, 9.0e-03)
      };
      \addlegendentry{low-order (staggering)}

      \addplot coordinates {
        (   980, 0.0797605  )
        (  3920, 0.00728726 )
        (  8820, 0.000574473)
        ( 15680, 2.99875e-05)
        ( 24500, 2.6495e-06 )
        ( 35280, 3.40046e-07)
        ( 48020, 6.96661e-08)
        ( 62720, 1.33605e-08)
        ( 79380, 4.99908e-09)
        ( 98000, 1.76148e-09)
        (118580, 4.59024e-10)
        (141120, 2.02394e-10)
        (165620, 8.03095e-11)
        (192080, 6.70143e-11)
      };
      \addlegendentry{high-order}
    \end{axis}
  \end{tikzpicture}}
\caption{\label{f:dbl-conv} The surfaces used in the Green's identity
  test in \pr{ss:result-dbl-layer} are shown at the top.  Below each
  surface, we plot the error $\|e\|_\infty$ as a function of the total
  number of grid points in the surface discretization.  The numerical
  results confirm that the low-order singularity subtraction scheme of
  Merkel converges as $\mathcal{O}(h^2)$ where
  $h = \mathcal{O}(1/\sqrt{N})$ is the mesh grid-spacing.  The
  high-order partition of unity based scheme is expected to converge
  exponentially and has a numerically observed convergence rate of
  roughly $12$th and $10$th order respectively for the two problems.
}
\end{figure}

\subsection{Normal component of the magnetic field due to external
  coils on the W7-X plasma surface using the virtual casing
  principle\label{ss:result-ext-field}}
We now present results for computing the normal component of the
external magnetic field using the numerical schemes based on
\pr{eq:BS_surface} and \pr{eq:vecpot_surface}.  In the setup shown in
\pr{f:W7X-Taylor}, we set $\Boundary$ to be the surface approximating
the plasma boundary in the W7-X stellarator and the magnetic field
$\vec{B}$ is a Taylor state with the surface~$\Boundary$ defined as
the last closed flux surface.  Taylor states are force-free magnetic
fields such that~$\nabla \times \vec{B} = \lambda \vec{B}$, where
$\lambda$ is a known constant called the Beltrami parameter.  We
construct this field using the method described
in~\cite{Malhotra2019}.  The reference solution
$\vec{B}_{\mathrm{ext}}\cdot\Normal$ is constructed by numerically
evaluating \pr{eq:BS_surface} on a very fine mesh ($\Nu=406$,
$\Nv=2030$) and we show convergence to this reference solution in
\pr{t:W7X-conv}.  In both formulations, we obtain spectral convergence
to nearly 10-digits.

  \begin{table}[b]
    \center
    \begin{tabular}{r @{\hskip 2em} c @{\hskip 2em} c}
      \toprule
      $\Nu~\times~\Nv$  & {Error with \pr{eq:BS_surface}} & {Error with \pr{eq:vecpot_surface}} \\
      \midrule
      $ 14~\times~~~70$ &   $1.6\nexp02$ &   $2.5\nexp02$ \\
      $ 42~\times~~210$ &   $4.2\nexp03$ &   $5.4\nexp04$ \\
      $ 84~\times~~420$ &   $2.3\nexp04$ &   $4.0\nexp06$ \\
      $112~\times~~560$ &   $4.5\nexp05$ &   $1.1\nexp06$ \\
      $154~\times~~770$ &   $7.1\nexp06$ &   $6.7\nexp08$ \\
      $196~\times~~980$ &   $2.4\nexp07$ &   $3.3\nexp09$ \\
      $252~\times~1260$ &   $5.3\nexp08$ &   $9.0\nexp10$ \\
      $322~\times~1610$ &   $4.0\nexp09$ &   $2.8\nexp10$ \\
      $392~\times~1960$ &   $2.0\nexp10$ &   $2.3\nexp10$ \\
      \bottomrule
    \end{tabular}
    \caption{\label{t:W7X-conv}
      Numerical results for the experimental setup shown in \pr{f:W7X-Taylor}
      using the two numerical schemes described by \pr{eq:BS_surface} and \pr{eq:vecpot_surface}.
      In each case, with increasing mesh resolution ($\Nu \times \Nv$), we report the relative error
      $\| \vec{B}^{0}_{ext} \cdot \Normal - \vec{B}_{ext} \cdot \Normal \|_{\infty} / \|\vec{B}\|_{\infty}$
      where $\vec{B}^{0}_{ext} \cdot \Normal$ is the reference solution.
      The reference solution is obtained numerically using
      \pr{eq:BS_surface} on a very fine mesh ($\Nu=406$, $\Nv=2030$).
    }
  \end{table}

  \begin{figure}[t]
    \center
    \includegraphics[width=0.75\textwidth]{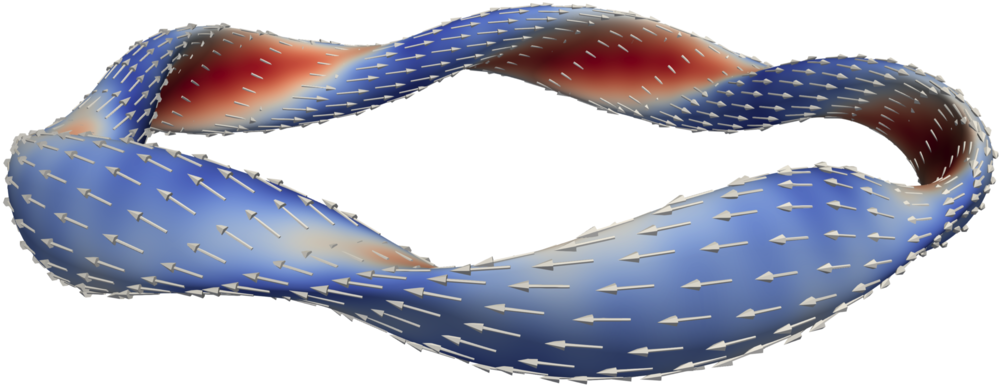}
    \caption{\label{f:W7X-Taylor}
      The plasma boundary $\Boundary$ in the W7-X stellarator is visualized here along with the magnetic field $\vec{B}$ for a fixed boundary calculation.
      The field $\vec{B}$ is a Taylor state with Beltrami parameter $\lambda=1$ and is generated using the method described in \cite{Malhotra2019}.
      Convergence results for computing the normal component of the external field for this setup are presented in \pr{t:W7X-conv}.
    }
  \end{figure}

\subsection{Vacuum field for the W7-X plasma surface and virtual casing principle}
  \label{subsec:results-vacuum-field}
  Consider the setup shown in \pr{f:W7X-loops}, in which the internal
  and external currents are given by known current loops $\gamma_1$
  and $\gamma_2$ respectively.  The magnetic field $\vec{B}$ produced
  by these currents at points $\vec{r}$ on $\Boundary$ is given by,
  \begin{equation}
    \vec{B}(\vec{r}) = \frac{\mu_{0}}{4\pi} \int_{\gamma_{1}} \frac{I_{1} \D\vec{l}' \times (\vec{r}-\vec{r}')}{|\vec{r}-\vec{r}'|^3}
                     + \frac{\mu_{0}}{4\pi} \int_{\gamma_{2}}
                     \frac{I_{2} \D\vec{l}' \times
                       (\vec{r}-\vec{r}')}{|\vec{r}-\vec{r}'|^3} ,
  \end{equation}
  where $I_1$ and $I_2$ are the currents in the loops $\gamma_1$ and
  $\gamma_2$ respectively.  Since the current loops do not intersect
  the surface $\Boundary$, the integrands are smooth and these
  integrals can be evaluated to high precision using standard
  quadratures.

  For our numerical test, we will now assume that we do not know the
  origin of the magnetic field thus calculated on the surface, and we
  wish to compute the normal component of the field due to the
  external loop~$\gamma_{2}$ using the virtual casing
  principle. Comparing our results using the virtual casing principle
  with the direct calculation using the Biot-Savart law would provide
  a reliable test of the accuracy of our quadrature scheme. This test
  would be complementary to the numerical test presented in
  Section~\ref{ss:result-ext-field}. However, we cannot directly apply
  the virtual casing principle\footnote{There is a more general version of
  the virtual casing
  principle\cite{drevlak2005pies,pustovitov2008decoupling,lazerson2012virtual,Hanson2015} 
  which can be applied in this case and would be more;
  however with this example we simply intend to show how the quadrature scheme
  discussed here can be used to solve boundary integral equations.}
  discussed in \pr{subsec:normal-fields}
  since it requires that $\Boundary$ be a flux surface
  (i.e. $\vec{B}\cdot\Normal = 0$ on $\Boundary$).  This issue can be
  easily circumvented by constructing a field
  $\vec{{B}}_{\mathrm{tot}} = \vec{B} + \nabla \Phi$ in
  $\Omega^{\mathrm{c}}$ such that
  $\vec{B}_{\mathrm{tot}} \cdot \Normal = 0$ on \Boundary. Physically,
  this is saying that one way of making $\Gamma$ a flux surface is to
  include an additional magnetic field source in~$\Omega$.  The
  corresponding potential~$\Phi$ can be constructed by solving
  \pr{eq:greensbasic} for the unknown $\Phi$ with given boundary conditions
  $\nabla \Phi \cdot \Normal$.
  This is a Fredholm integral equation of the second kind
  and can therefore be solved efficiently using GMRES as described in
  \pr{subsec:vacuum-fields}.
  The numerical test in this section is
  therefore a combined test of our quadrature scheme for the two
  applications discussed in Section \ref{sec:two_problems}: the
  virtual casing principle and the calculation of vacuum fields.  The
  normal component
  data~$\nabla \Phi \cdot \Normal = -\vec{B} \cdot \Normal$ is known
  a priori
  from the construction of $\Phi$, and the tangential components of
  $\nabla \Phi$ are computed using Fourier differentiation of $\Phi$
  on the surface.  The virtual casing principle can now be be applied
  to $\vec{{B}}_{\mathrm{tot}}$.
  We use \pr{eq:BS_surface} and \pr{eq:vecpot_surface} to compute
  $\vec{B}_{\mathrm{ext}} \cdot \Normal$, where
  $\vec{B}_{\mathrm{ext}}$ is the field due to the external currents,
  namely the external current loop $\gamma_{2}$ in this example. We
  compare this expression with the reference normal component of the
  field (computed directly via Biot-Savart applied to the external
  current loop~$\gamma_2$) and report convergence results in
  \pr{t:vacuum-conv}. We observe spectral convergence as we refine the
  mesh and reduce the tolerance $\gmrestol$ for the GMRES solve.  In
  addition, $\gmresiter$ remains reasonably small since we are solving
  second-kind integral equations.
  
  Finally, we observe that with a generalized version of the virtual casing principle \cite{Hanson2015}, it is possible to compute all the components of the magnetic field on $\Gamma$ due to the external current loop $\gamma_{2}$ with only the magnetic field $\mathbf{B}$ on $\Gamma$ given as input. We implemented this more general numerical test as well, in which we compare all the components of $\vec{B}_{\mathrm{ext}}$, and obtained similar results to the ones reported for the normal component here. We chose to focus on the normal component in this section because it is the most common application in magnetic confinement fusion.

  \begin{figure}[t]
    \center
    \resizebox{.80\linewidth}{!}{\begin{picture}(300,110)
      \put(0,0){\includegraphics[width=300pt]{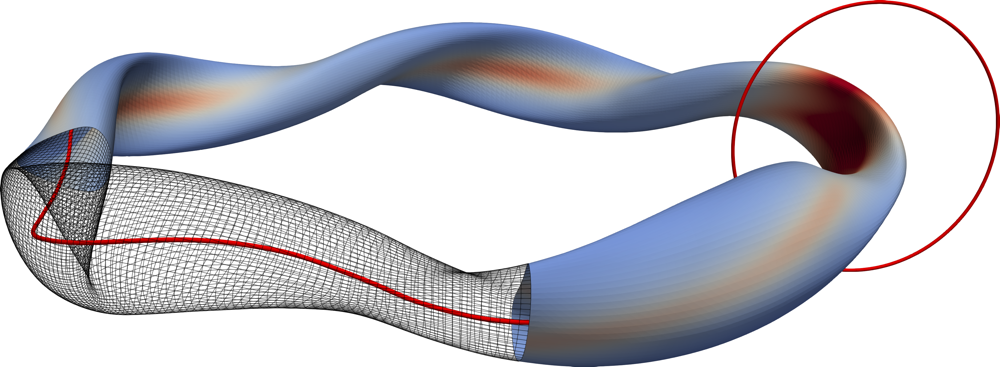}}
      \put(60,25){\color{red} \huge $\gamma_{1}$}
      \put(265,95){\color{red} \huge $\gamma_{2}$}
    \end{picture}}
    \caption{\label{f:W7X-loops}
      The boundary $\Boundary$ along with the internal ($\gamma_1$) and external ($\gamma_1$) current loops.
      The current loops produce a magnetic field \vec{B} and the magnitude of this field on $\Boundary$ is shown.
      Using the methods described in \pr{subsec:normal-fields,subsec:vacuum-fields}, given only the magnetic field \vec{B} on $\Boundary$, we can recover the magnetic field on $\Boundary$ due to the external current loop $\gamma_2$.
      In \pr{t:vacuum-conv}, we show convergence results by comparing the external field recovered from \vec{B} with the field computed directly from integrating over the external current loop $\gamma_2$.
    }
  \end{figure}

  \begin{table}[b]
    \center
    \begin{tabular}{r @{\hskip 2em} r @{\hskip 2em} r @{\hskip 2em} c @{\hskip 2em} c}
      \toprule
      $\Nu~\times~\Nv$  & $\gmrestol$ & $\gmresiter$ & {Error with \pr{eq:BS_surface}} & {Error with \pr{eq:vecpot_surface}} \\
      \midrule
      $ 14~\times~~~70$ &  $1\nexp01$ &         $ 3$ &                    $1.2\nexp01$ &                        $1.0\nexp01$ \\
      $ 42~\times~~210$ &  $3\nexp03$ &         $ 6$ &                    $2.2\nexp02$ &                        $1.5\nexp02$ \\
      $ 84~\times~~420$ &  $1\nexp04$ &         $ 8$ &                    $9.2\nexp04$ &                        $2.8\nexp04$ \\
      $112~\times~~560$ &  $3\nexp05$ &         $ 9$ &                    $1.0\nexp04$ &                        $5.2\nexp05$ \\
      $154~\times~~770$ &  $1\nexp06$ &         $13$ &                    $2.0\nexp05$ &                        $5.6\nexp06$ \\
      $196~\times~~980$ &  $1\nexp07$ &         $14$ &                    $1.1\nexp06$ &                        $5.3\nexp07$ \\
      $252~\times~1260$ &  $1\nexp08$ &         $16$ &                    $9.6\nexp08$ &                        $7.3\nexp08$ \\
      $322~\times~1610$ &  $1\nexp09$ &         $18$ &                    $1.6\nexp08$ &                        $1.1\nexp08$ \\
      $392~\times~1960$ &  $1\nexp10$ &         $20$ &                    $1.2\nexp09$ &                        $8.3\nexp10$ \\
      \bottomrule
    \end{tabular}
    \caption{\label{t:vacuum-conv}
      Numerical results for the experimental setup described in \pr{subsec:results-vacuum-field}.
      As we refine the mesh and reduce the tolerance $\gmrestol$ for the GMRES solve,
      we observe spectral convergence in the relative error
      $\| \vec{B}^{0}_{ext} \cdot \Normal - \vec{B}_{ext} \cdot \Normal \|_{\infty} / \|\vec{B}\|_{\infty}$
      where $\vec{B}^{0}_{ext} \cdot \Normal$ is the reference solution computed by directly integrating along the external current loop $\gamma_2$.
      For each case, we also report the number of GMRES iterations $\gmresiter$ in computing the potential $\phi$.
    }
  \end{table}

\section{Summary}
\label{sec:summary}

We have presented a numerical quadrature scheme for the evaluation of
integrals with singular kernels appearing in magnetostatic problems
for magnetic confinement fusion applications. We demonstrated that our
scheme has high-order convergence, and is much more accurate than the
schemes currently used in the magnetic fusion community, even for the
coarsest meshes. Implementing our algorithm in current magnetic fusion
codes based on integral formulations with singular kernels would
significantly improve their accuracy for a fixed mesh size;
conversely, for a fixed target accuracy, our scheme would allow these
codes to use significantly coarser meshes. Although we have not done
so for the numerical tests presented in this article, our scheme can
be coupled with fast multipole methods to obtain optimal asymptotic
scaling for sufficiently large-scale problems.

There are two immediate extensions to the work presented in this
article. The first is the development of a fast, high-order quadrature
scheme for the evaluation of integrals with singular kernels on
surfaces with edges, which is needed for situations in which the
surface has one or several magnetic separatrices. The second extension
is the construction of a fast and high-order quadrature scheme for the
evaluation of integrals with singular kernels at points away from the
surface but close to it (there is no particular numerical issue for
points far away from the surface). Such a scheme would, for example,
be necessary if we wanted to evaluate the off-surface vacuum field
potential~$\Phi$ given by (\ref{eq:greensoutside}), or to compute the
off-surface magnetic field near a bounding surface in the integral
equation formulation of the calculation of Taylor states as
implemented in BIEST \cite{Malhotra2019}. Both extensions reflect
questions which have not yet been addressed in a definitive way by the
applied mathematics community~\cite{Helsing_2008, Helsing_2011,
  Helsing_2013, Serkh_2016, Serkh_2019, Hoskins_2019, Zhao_2010,
  Corona_2017, Bremer_2012, Bremer_2013, Tlupova_2013, Kl_ckner_2013,
  Wala_2019, Af_Klinteberg_2014, Siegel_2018, Ying_2006,
  Rahimian_2017, Carvalho_2018a, Carvalho_2018b}.
These problems are the subject of ongoing work, with results to be
reported in the future.

\section{Acknowledgements}
The authors would like to thank Allen Boozer for bringing to our
attention the
question treated in this article. D.M. was partially supported by the Office of Naval Research under award number \#N00014-17-1-2451, the Simons Foundation/SFARI (560651, AB), and the U.S. Department of Energy, Office of Science, Fusion Energy Sciences under Awards No. DE-FG02-86ER53223 and DE-SC0012398. A.J.C. was partially supported by the Simons Foundation/SFARI (560651, AB), and the U.S. Department of Energy, Office of Science, Fusion Energy Sciences under Award Nos. DE-FG02-86ER53223 and DE-SC0012398. M.O. was partially supported by the Office of Naval Research under award numbers \#N00014-17-1-2451 and \#N00014-17-1-2059, and the Simons Foundation/SFARI (560651, AB). E.T. was partially supported by the U.S. Department of Energy, Office of Science, Fusion Energy Sciences under Award Nos. DE-FG02-86ER53223 and DE-SC0012398.

\bibliographystyle{abbrv}
\bibliography{ref}

\begin{thebibliography}{10}

\bibitem{Af_Klinteberg_2014}
L.~af~Klinteberg and A.-K. Tornberg.
\newblock Fast {Ewald} summation for {Stokesian} particle suspensions.
\newblock {\em International Journal for Numerical Methods in Fluids},
  76(10):669--698, Sep 2014.

\bibitem{Albanese1988}
R.~Albanese and G.~Rubinacci.
\newblock Integral formulation for {3D} eddy-current computation using edge
  elements.
\newblock {\em IEE Proceedings A - Physical Science, Measurement and
  Instrumentation, Management and Education - Reviews}, 135(7):457--462, Sep
  1988.

\bibitem{Atanasiu1999}
C.~V. Atanasiu, A.~H. Boozer, L.~E. Zakharov, A.~A. Subbotin, and G.~I. Miron.
\newblock Determination of the vacuum field resulting from the perturbation of
  a toroidally symmetric plasma.
\newblock {\em Physics of Plasmas}, 6(7):2781--2790, Jul 1999.

\bibitem{Bremer_2012}
J.~Bremer and Z.~Gimbutas.
\newblock A {Nyström} method for weakly singular integral operators on
  surfaces.
\newblock {\em Journal of Computational Physics}, 231(14):4885--4903, May 2012.

\bibitem{Bremer_2013}
J.~Bremer and Z.~Gimbutas.
\newblock On the numerical evaluation of the singular integrals of scattering
  theory.
\newblock {\em Journal of Computational Physics}, 251:327--343, Oct 2013.

\bibitem{Bruno_2001a}
O.~P. Bruno and L.~A. Kunyansky.
\newblock A fast, high-order algorithm for the solution of surface scattering
  problems: Basic implementation, tests, and applications.
\newblock {\em Journal of Computational Physics}, 169(1):80--110, May 2001.

\bibitem{Bruno_2001b}
O.~P. Bruno and L.~A. Kunyansky.
\newblock Surface scattering in three dimensions: an accelerated high-order
  solver.
\newblock {\em Proceedings of the Royal Society of London. Series A:
  Mathematical, Physical and Engineering Sciences}, 457(2016):2921--2934, Dec
  2001.

\bibitem{Carvalho_2018b}
C.~{Carvalho}, S.~{Khatri}, and A.~D. {Kim}.
\newblock {Asymptotic approximations for the close evaluation of double-layer
  potentials}.
\newblock {\em arXiv e-prints}, page arXiv:1810.02483, Oct 2018.

\bibitem{Carvalho_2018a}
C.~{Carvalho}, S.~{Khatri}, and A.~D. {Kim}.
\newblock {Close evaluation of layer potentials in three dimensions}.
\newblock {\em arXiv e-prints}, page arXiv:1807.02474, Jul 2018.

\bibitem{Chance2007}
M.~Chance, A.~Turnbull, and P.~Snyder.
\newblock Calculation of the vacuum {Green}’s function valid even for high
  toroidal mode numbers in tokamaks.
\newblock {\em Journal of Computational Physics}, 221(1):330--348, Jan 2007.

\bibitem{Chance1997}
M.~S. Chance.
\newblock Vacuum calculations in azimuthally symmetric geometry.
\newblock {\em Physics of Plasmas}, 4(6):2161--2180, Jun 1997.

\bibitem{Corona_2017}
E.~Corona, L.~Greengard, M.~Rachh, and S.~Veerapaneni.
\newblock An integral equation formulation for rigid bodies in {Stokes} flow in
  three dimensions.
\newblock {\em Journal of Computational Physics}, 332:504--519, Mar 2017.

\bibitem{Drevlak2018}
M.~Drevlak, C.~Beidler, J.~Geiger, P.~Helander, and Y.~Turkin.
\newblock Optimisation of stellarator equilibria with {ROSE}.
\newblock {\em Nuclear Fusion}, 59(1):016010, Nov 2018.

\bibitem{drevlak2005pies}
M.~Drevlak, D.~Monticello, and A.~Reiman.
\newblock Pies free boundary stellarator equilibria with improved initial
  conditions.
\newblock {\em Nuclear fusion}, 45(7):731, 2005.

\bibitem{Freidberg1975}
J.~P. Freidberg and W.~Grossmann.
\newblock Magnetohydrodynamic stability of a sharp boundary model of tokamak.
\newblock {\em The Physics of Fluids}, 18(11):1494--1506, 1975.

\bibitem{Freidberg1976}
J.~P. Freidberg, W.~Grossmann, and F.~A. Haas.
\newblock Stability of a high-{$\beta$}, l=3 stellarator.
\newblock {\em The Physics of Fluids}, 19(10):1599--1607, 1976.

\bibitem{Gruber19812}
R.~Gruber, S.~Semenzato, F.~Troyon, T.~Tsunematsu, W.~Kerner, P.~Merkel, and
  W.~Schneider.
\newblock {Hera} and other extensions of {Erato}.
\newblock {\em Computer Physics Communications}, 24(3):363--376, Dec 1981.

\bibitem{Gruber1981}
R.~Gruber, F.~Troyon, D.~Berger, L.~Bernard, S.~Rousset, R.~Schreiber,
  W.~Kerner, W.~Schneider, and K.~Roberts.
\newblock {Erato} stability code.
\newblock {\em Computer Physics Communications}, 21(3):323--371, Jan 1981.

\bibitem{guentherlee1995}
R.~B. Guenther and J.~W. Lee.
\newblock {\em Partial Differential Equations of Mathematical Physics and
  Integral Equations}.
\newblock New York: Dover Publications, Inc., 1995.

\bibitem{Hanson2015}
J.~D. Hanson.
\newblock The virtual-casing principle and {Helmholtz}'s theorem.
\newblock {\em Plasma Physics and Controlled Fusion}, 57(11):115006, Sep 2015.

\bibitem{Helsing_2011}
J.~Helsing.
\newblock A fast and stable solver for singular integral equations on piecewise
  smooth curves.
\newblock {\em SIAM Journal on Scientific Computing}, 33(1):153--174, Jan 2011.

\bibitem{helsing2015explicit}
J.~Helsing and A.~Karlsson.
\newblock An explicit kernel-split panel-based {Nyström} scheme for integral
  equations on axially symmetric surfaces.
\newblock {\em Journal of Computational Physics}, 272:686--703, Sep 2014.

\bibitem{Helsing_2008}
J.~Helsing and R.~Ojala.
\newblock Corner singularities for elliptic problems: Integral equations,
  graded meshes, quadrature, and compressed inverse preconditioning.
\newblock {\em Journal of Computational Physics}, 227(20):8820--8840, Oct 2008.

\bibitem{Helsing_2013}
J.~Helsing and K.-M. Perfekt.
\newblock On the polarizability and capacitance of the cube.
\newblock {\em Applied and Computational Harmonic Analysis}, 34(3):445--468,
  May 2013.

\bibitem{Hirshman1986}
S.~Hirshman, W.~van RIJ, and P.~Merkel.
\newblock Three-dimensional free boundary calculations using a spectral
  {Green}'s function method.
\newblock {\em Computer Physics Communications}, 43(1):143--155, Dec 1986.

\bibitem{Hoskins_2019}
J.~G. Hoskins, V.~Rokhlin, and K.~Serkh.
\newblock On the numerical solution of elliptic partial differential equations
  on polygonal domains.
\newblock {\em SIAM Journal on Scientific Computing}, 41(4):A2552--A2578, Jan
  2019.

\bibitem{Kl_ckner_2013}
A.~Kl{\"o}ckner, A.~Barnett, L.~Greengard, and M.~O'Neil.
\newblock Quadrature by expansion: A new method for the evaluation of layer
  potentials.
\newblock {\em Journal of Computational Physics}, 252:332--349, Nov 2013.

\bibitem{Lackner1976}
K.~Lackner.
\newblock Computation of ideal {MHD} equilibria.
\newblock {\em Computer Physics Communications}, 12(1):33--44, Sep 1976.

\bibitem{Landreman2017}
M.~Landreman.
\newblock An improved current potential method for fast computation of
  stellarator coil shapes.
\newblock {\em Nuclear Fusion}, 57(4):046003, Feb 2017.

\bibitem{Lazanja2011}
D.~Lazanja.
\newblock {\em Components of the magnetic field from interior and exterior
  sources in toroidal geometry}.
\newblock PhD thesis, Columbia University, 2011.

\bibitem{lazerson2012virtual}
S.~A. Lazerson.
\newblock The virtual-casing principle for 3d toroidal systems.
\newblock {\em Plasma Physics and Controlled Fusion}, 54(12):122002, 2012.

\bibitem{Lee2015}
J.~P. Lee, A.~Cerfon, J.~P. Freidberg, and M.~Greenwald.
\newblock Tokamak elongation - how much is too much? part 2. numerical
  results.
\newblock {\em Journal of Plasma Physics}, 81(6):515810608, Dec 2015.

\bibitem{Malhotra2019}
D.~Malhotra, A.~Cerfon, L.-M. Imbert-Gérard, and M.~O'Neil.
\newblock {Taylor} states in stellarators: A fast high-order boundary integral
  solver.
\newblock {\em Journal of Computational Physics}, 397:108791, Nov 2019.

\bibitem{Marx2017}
A.~Marx and H.~Lütjens.
\newblock Free-boundary simulations with the {XTOR}-2f code.
\newblock {\em Plasma Physics and Controlled Fusion}, 59(6):064009, May 2017.

\bibitem{Merkel1982}
P.~Merkel.
\newblock A {Green}'s function method for the vacuum contribution to the {MHD}
  stability of helically symmetric equilibria.
\newblock {\em Zeitschrift fur Naturforschung - Section A Journal of Physical
  Sciences}, 37(8):859--865, Aug 1982.

\bibitem{Merkel1986}
P.~Merkel.
\newblock An integral equation technique for the exterior and interior
  {Neumann} problem in toroidal regions.
\newblock {\em Journal of Computational Physics}, 66(1):83--98, Sep 1986.

\bibitem{Merkel1987}
P.~Merkel.
\newblock Solution of stellarator boundary value problems with external
  currents.
\newblock {\em Nuclear Fusion}, 27(5):867--871, May 1987.

\bibitem{Morozov1966}
A.~I. {Morozov} and L.~S. {Solov'ev}.
\newblock {Motion of Charged Particles in Electromagnetic Fields}.
\newblock {\em Reviews of Plasma Physics}, 2:201, Jan 1966.

\bibitem{oneil2018taylor}
M.~O'Neil and A.~J. Cerfon.
\newblock An integral equation-based numerical solver for {Taylor} states in
  toroidal geometries.
\newblock {\em Journal of Computational Physics}, 359:263--282, Apr 2018.

\bibitem{pustovitov2001magnetic}
V.~Pustovitov.
\newblock Magnetic diagnostics: General principles and the problem of
  reconstruction of plasma current and pressure profiles in toroidal systems.
\newblock {\em Nuclear fusion}, 41(6):721, 2001.

\bibitem{pustovitov2008decoupling}
V.~Pustovitov.
\newblock Decoupling in the problem of tokamak plasma response to asymmetric
  magnetic perturbations.
\newblock {\em Plasma Physics and Controlled Fusion}, 50(10):105001, 2008.

\bibitem{Pustovitov2008}
V.~D. Pustovitov.
\newblock General formulation of the resistive wall mode coupling equations.
\newblock {\em Physics of Plasmas}, 15(7):072501, Jul 2008.

\bibitem{Rachh_2017}
M.~Rachh, A.~Klöckner, and M.~O~'Neil.
\newblock Fast algorithms for {Quadrature by Expansion} {I}: Globally valid
  expansions.
\newblock {\em Journal of Computational Physics}, 345:706--731, Sep 2017.

\bibitem{Rahimian_2017}
A.~Rahimian, A.~Barnett, and D.~Zorin.
\newblock Ubiquitous evaluation of layer potentials using {Quadrature by
  Kernel-Independent Expansion}.
\newblock {\em BIT Numerical Mathematics}, 58(2):423--456, Nov 2017.

\bibitem{Serkh_2019}
K.~Serkh.
\newblock On the solution of elliptic partial differential equations on regions
  with corners {II}: Detailed analysis.
\newblock {\em Applied and Computational Harmonic Analysis}, 46(2):250--287,
  Mar 2019.

\bibitem{Serkh_2016}
K.~Serkh and V.~Rokhlin.
\newblock On the solution of elliptic partial differential equations on regions
  with corners.
\newblock {\em Journal of Computational Physics}, 305:150--171, Jan 2016.

\bibitem{Shafranov1972}
V.~Shafranov and L.~Zakharov.
\newblock Use of the virtual-casing principle in calculating the containing
  magnetic field in toroidal plasma systems.
\newblock {\em Nuclear Fusion}, 12(5):599--601, Sep 1972.

\bibitem{Siegel_2018}
M.~Siegel and A.-K. Tornberg.
\newblock A local target specific quadrature by expansion method for evaluation
  of layer potentials in {3D}.
\newblock {\em Journal of Computational Physics}, 364:365--392, Jul 2018.

\bibitem{stratton1941electromagnetic}
J.~A. Stratton.
\newblock Electromagnetic theory, mcgrow-hill book company.
\newblock {\em Inc., New York, and London}, 1941.

\bibitem{Strumberger1997}
E.~Strumberger.
\newblock Finite-$\beta$ magnetic field line tracing for {Helias}
  configurations.
\newblock {\em Nuclear Fusion}, 37(1):19--27, Jan 1997.

\bibitem{Tlupova_2013}
S.~Tlupova and J.~T. Beale.
\newblock Nearly singular integrals in {3D} {Stokes} flow.
\newblock {\em Communications in Computational Physics}, 14(5):1207--1227, Nov
  2013.

\bibitem{Trefethen2014}
L.~N. Trefethen and J.~A.~C. Weideman.
\newblock The exponentially convergent trapezoidal rule.
\newblock {\em SIAM Review}, 56(3):385--458, Jan 2014.

\bibitem{Wala_2019}
M.~Wala and A.~Klöckner.
\newblock A fast algorithm for {Quadrature by Expansion} in three dimensions.
\newblock {\em Journal of Computational Physics}, 388:655--689, Jul 2019.

\bibitem{Ying_2006}
L.~Ying, G.~Biros, and D.~Zorin.
\newblock A high-order {3D} boundary integral equation solver for elliptic pdes
  in smooth domains.
\newblock {\em Journal of Computational Physics}, 219(1):247--275, Nov 2006.

\bibitem{Zhao_2010}
H.~Zhao, A.~H. Isfahani, L.~N. Olson, and J.~B. Freund.
\newblock A spectral boundary integral method for flowing blood cells.
\newblock {\em Journal of Computational Physics}, 229(10):3726--3744, May 2010.

\end{thebibliography}
\end{document}